\documentclass[twoside]{article}
\date{}
\usepackage{graphicx}

\usepackage{amssymb,amsfonts}

\def\NN{\mathbb{N}}    
\def\ZZ{\mathbb{Z}}    
\def\PP{\mathbb{P}}    
\def\EE{\mathbb{E}}    

\def\proof{\noindent{\em Proof.}~}
\def\endproof{\hfill\mbox{\fbox{\rule{0mm}{0.0mm}}}\vspace{2ex}}
\newtheorem{theorem}{Theorem}
\newtheorem{lemma}[theorem]{Lemma}
\newtheorem{proposition}[theorem]{Proposition}
\newtheorem{corollary}[theorem]{Corollary}
\title{A sharp uniform bound for the distribution of sums of Bernoulli trials}
\author{
Jean-Bernard Baillon\\[-0.5ex]
{\small SAMM--EA 4543, Universit\'e de Paris 1}\\[-0.5ex]
{\small 75013 Paris, France }\\[1ex]
\and
Roberto Cominetti\\[-0.5ex] 
{\small Departamento de Ingenier\'\i a Industrial, Universidad de Chile}\\[-0.5ex] 
{\small Avda. Rep\'ublica 701, Santiago, Chile}\\[1ex]
  \and 
Jos\'e Vaisman\\[-0.5ex] 
{\small Departamento de Ingenier\'\i a Matem\'atica, Universidad de Chile}\\[-0.5ex]
{\small Avda. Blanco Encalada 2120, Santiago, Chile}
}
\begin{document}
\maketitle

\begin{abstract}
In this note we establish a uniform bound for the distribution of a sum
$S_n=X_1+\cdots+X_n$ of independent non-homogeneous Bernoulli trials. 
Specifically, we prove that 
$\sigma_n\,\PP(S_n\!=\!j)\!\leq\! \eta$ where 
$\sigma_n$ denotes the standard deviation of $S_n$ and $\eta$ is a universal constant.
We compute the best possible constant $\eta\!\sim\! 0.4688$ and we show that
the bound also holds for limits of sums and differences of Bernoullis, including the Poisson 
laws which constitute the worst case and attain the bound. We also investigate the
optimal bounds for $n$ and $j$ fixed.
An application to estimate the rate of convergence of Mann's fixed point iterations 
is presented.
\end{abstract}

\vspace{2ex}


\vspace{2ex}
{\small \noindent{\bf Running title:}  A sharp uniform bound for sums of Bernoullis} 
\vfill

\section{Introduction}
Let $S_n=X_1+\cdots+X_n$ be a sum of independent non-homogeneous Bernoulli trials with success 
probabilities $p_i$. The distribution of $S_n$ is known to be unimodal and bell-shaped with mean 
$\mu_n=\sum_{i=1}^np_i$ and variance $\sigma_n^2=\sum_{i=1}^np_i(1\!-\!p_i)$.
Its mode is either $\lfloor\mu_n\rfloor$ or $\lceil\mu_n\rceil$ or both \cite{dar,sam}, and the same
holds for the median \cite{jog}. In this paper we investigate how large it can be the modal probability.
More precisely, we establish a uniform upper bound
\begin{equation}\label{bound}
\sigma_n\,\PP(S_n\!\!=\!j)\leq \eta
\end{equation}
for all $n$, $j$ and $p_i$, and we prove that the best possible constant is
\begin{equation}\label{bnd1}
\eta=\max_{\lambda\geq 0}\sqrt{2\lambda}\;e^{-2\lambda}\sum_{k=0}^\infty\left({\lambda^k\over k!}\right)^2\sim0.4688.
\end{equation}

The existence of a universal bound (\ref{bound}) can be established using tools related to the local 
limit theorem \cite{bmd,gam,mcdo}. It also follows as a special case of the Kolmogorov-Rogozin concentration 
inequality \cite{rog} which states  a more general bound valid for discrete random variables $X_i$ with 
$\sigma_n$ replaced by $\sqrt{\sum_{i=1}^n(1\!-\!\psi_i)}$ where $\psi_i=\max_x\PP(X_i=x)$. For sums of 
Bernoullis this is equivalent to (\ref{bound}), so that our contribution is mainly the computation of the 
optimal constant $\eta$, as well as the identification 
of the role of the Poisson law in the worst case.
Namely, the local limit theorem shows that for a wide range of random variables the limit as $n\to \infty$ in (\ref{bound})  exists and equals $1/\sqrt{2\pi}$. 
Since $\eta$ exceeds this value, it follows that  the worst case situation is not associated with random 
variables obeying the central limit theorem. It is then natural to expect that the worst case may have 
to do with the Poisson law, and that is what actually happens. 
In fact the expression (\ref{bnd1}) is just 
\begin{equation}\label{bnd2}
\eta=\max_{\lambda\geq 0}\sqrt{2\lambda}\;\PP(N_\lambda\!=\!N'_\lambda)
\end{equation}
where $N_\lambda$ and $N'_\lambda$ are independent Poisson variables with parameter $\lambda$.
   Since the Poisson law also happens to be extremal in other 
bounds such as Rosenthal's inequality ({\em cf.} \cite{fig,ibr,ssh,ute}), a natural question is whether (\ref{bound})  
might hold for more general sums of random variables.

The inequality (\ref{bound}) complements the large deviation bounds
that provide estimates of the form $\PP(|S_n\!\!-\EE(S_n)|\geq t)\leq f(nt^2)$
with $f(x)\to 0$ as $x\to\infty$, usually at an exponential rate
 ({\em cf.} \cite{alo,azu,ber,che,fel,goh,hoe,mcd,pet}). In contrast, (\ref{bound}) does not 
give such fast asymptotic rates but it can be used to bound 
$\PP(S_n\!\!=\!j)$ for all values of $j$ including values close to the mean $\EE(S_n)$.
This  already proved useful in establishing an optimality 
guarantee for an approximation algorithm in discrete stochastic optimization
(see \cite{ccs}). In this paper we present another application to 
the rate of convergence of fixed point iterations for non-expansive maps.
In both settings a sharp constant $\eta$ is relevant as it yields better bounds.

The paper is organized as follows. In \S\ref{sharp} we show the sharp 
uniform bound (\ref{bound}) to be valid for all $n, j$ and $p_i$'s, and we briefly
discuss some extensions to more general distributions including 
sums and differences of Bernoullis as well as their limits which cover all 
Poisson distributions and more.  
In \S\ref{s2} we investigate more closely the optimal bounds for fixed $n$ and $j$.
In the final section \S \ref{appl} we show how (\ref{bound}) allows to establish 
the rate of convergence for fixed point iterations.

\section{A sharp uniform bound}\label{sharp}
 
 \begin{theorem}\label{T}
Let $S_n\!=\!X_1+\cdots+X_n$ be a sum of independent Bernoulli
trials with $\PP(X_i\!=\!1)=p_i$, and let $\sigma_n^2\!=\!\sum_{i=1}^np_i(1\!-\!p_i)$
denote its variance. Then
\begin{equation}\label{bound1}
\sigma_n\;\PP(S_n\!\!=\!j)\leq \eta
\end{equation}
where $\eta=\max_{\lambda>0}\sqrt{2\lambda}\;\PP(N_\lambda\!=\!N'_\lambda)$ with 
$N_\lambda$ and $N'_\lambda$ independent Poisson variables of parameter $\lambda$.
This bound is sharp and we have more explicitly $\eta=\max_{x\geq 0}\sqrt{x}\,e^{-x}I_0(x)\sim 0.4688$,
where $I_0(x)$ is the modified Bessel function
$$I_0(x)=\mbox{$\sum_{k=0}^\infty({x^k\over 2^kk!})^2$}=\frac{1}{\pi}\int_0^{\pi}\!\!\exp(x\cos\theta)\,d\theta.$$ 
\end{theorem}

\proof Consider the generating function
$\phi(z)=\EE(z^{S_n})=\sum_{j=0}^n\PP(S_n\!=\!j)z^j$. Integrating $\phi(z)/z^{j+1}$ along the 
unit circle $\mathcal{C}$ in the complex plane we get
\begin{equation}\label{pgf}
\PP(S_n\!=\!j)=\mbox{${1\over 2\pi i}\int_{\mathcal{C}}{\phi(z)\over z^{j+1}}\,dz=
{1\over 2\pi}\int_0^{2\pi}\!\phi(e^{i\theta})e^{-ij\theta}\,d\theta$}
\end{equation}
so that taking absolute value it follows that
\begin{equation}\label{pgf2}
\PP(S_n\!=\!j)\leq\mbox{${1\over 2\pi}\int_0^{2\pi}\!|\phi(e^{i\theta})|\,d\theta$}.
\end{equation}
The independence of the $X_i$'s yields $\phi(z)=\EE[\prod_{i=1}^nz^{X_i}]=\prod_{i=1}^n[(1\!-\!p_i)+p_i z]$,
from which we obtain
$$|\phi(e^{i\theta})|=\prod_{i=1}^n\sqrt{1+2p_i(1\!-\!p_i)(\cos\theta\!-\!1)}.$$
Using the inequality $1+y_i\leq\exp(y_i)$ with $y_i=2p_i(1\!-\!p_i)(\cos\theta\!-\!1)$,
and setting $x=\sum_{i=1}^np_i(1\!-\!p_i)$, we deduce
$$|\phi(e^{i\theta})|\leq\prod_{i=1}^n\exp(y_i/2)=\exp(x(\cos\theta\!-\!1)),$$
and since $\sigma_n=\sqrt{x}$ we conclude
$$\sigma_n\,\PP(S_n=j)\leq
\sqrt{x}\mbox{$\frac{1}{2\pi}\int_0^{2\pi}\exp(x(\cos\theta\!-\!1))\,d\theta$}=\sqrt{x}\,e^{-x}I_0(x)\leq\eta.$$

In order to show that the bound is sharp, consider a sum of $n=2a$ Bernoullis, half of them with 
$p_i=\lambda/a$ and the other half with $p_i=1\!-\!\lambda/a$, so that $S_n=U+V$ with
$U\sim B(a,\lambda/a)$ and $V\sim B(a,1\!-\!\lambda/a)$ independent Binomials.
Note that $V\stackrel{d}{=}a-U'$ with $U'$ an independent copy of $U$,
and for $j=a$ we get
$$\sigma_n\,\PP(S_n\!=\!a)=
\mbox{$\sqrt{2\lambda(1\!-\!{\lambda\over a})}\;\PP(U=U').$}$$
Since $U$ and $U'$ converge as $a\to\infty$ to independent Poisson variables $N_{\lambda}$ 
and $N'_{\lambda}$ with parameter $\lambda$, this expression tends to
$$\sqrt{2\lambda}\;\PP(N_{\lambda}=N'_{\lambda})=\sqrt{2\lambda}\,e^{-2\lambda}\sum_{k=0}^\infty(\mbox{${\lambda^k\over k!}$})^2= \sqrt{2\lambda}\,e^{-2\lambda}I_0(2\lambda)$$
which proves that the bound is sharp.
\endproof

\noindent{\sc Remark 1.}
The optimal bound $\eta=\max_{\lambda>0}\sqrt{2\lambda}\;\PP(N_\lambda\!=\!N'_\lambda)$
is approximately $\eta\sim 0.468822355499$ which is attained for $\lambda\sim 0.39498893$.

\vspace{2ex}
\noindent{\sc Remark 2.} The proof above shows that the bound $\eta$ is 
asymptotically attained for a sum of two Binomials with different success probabilities
$p=\lambda/a$ and $p'=1-\lambda/a$. As a matter of fact, allowing for two different 
Binomials is essential since for a single Binomial $S_n\sim B(n,x)$ we have the sharper bound
\begin{equation}\label{single}
\sigma_n\,\PP(S_n\!=\!j)=\mbox{$\sqrt{n\,x(1\!-\!x)}\;{n\choose j}x^j(1\!-\!x)^{n-j}\leq{1\over \sqrt{2e}}$}
\end{equation}
with ${1\over \sqrt{2e}}\sim 0.4289 <\eta$.
To prove (\ref{single}) we note that for $n$ and $j$ given, the maximum over $x\in [0,1]$
 is attained at
$x=(j+{1\over 2})/(n+1)$, so that replacing this value all we must show is  that $C^n_j\leq {1\over\sqrt{2e}}$ where
$$C^n_j={n\choose j}{\sqrt{n}(j+\mbox{${1\over 2}$})^{j+{1\over 2}}(n-j+\mbox{${1\over 2}$})^{n-j+{1\over 2}}\over (n+1)^{n+1}}.$$
Now $C^n_{j+1}/C^n_j=H(n-j)/H(j+1)$ where 
$H(x)=x(x-\mbox{${1\over 2}$})^{x-\mbox{${1\over 2}$}}/(x+\mbox{${1\over 2}$})^{x+\mbox{${1\over 2}$}}$
is decreasing, so
that $C^n_j$ decreases for $j\leq {n-1\over 2}$ and increases 
afterwards. Hence $C_j^n$ is maximal at $j=0$ or $j=n$, and then
the conclusion follows since
$$C^n_0=C^n_n=
{1\over\sqrt{2}}\sqrt{n\over n+{1\over 2}}\left(1-{1\over 2(n+1)}\right)^{n+1}\leq {1\over\sqrt{2}}\exp(-\mbox{${1\over 2}$})={1\over\sqrt{2e}}.
$$

\subsection{Extension to more general distributions}\label{extension}
As a consequence of Theorem 1 we see that (\ref{bound}) still holds 
for any random variable $S_n\!=\!\sum_{i=1}^n\pm X_i$ that can be expressed as {\em sums and differences} 
of independent Bernoullis. Moreover, the bound remains true for
limits of such variables, which includes all Poisson distributions
as well as infinite series
$S^\infty=\sum_{i=1}^\infty X_i$ of independent Bernoullis
with $\sum_{i=1}^\infty p_i<\infty$, namely


\begin{corollary}\label{cor2} Let $S=(X+S^\infty_+)-(Y+S^\infty_-)$ with $X,Y$ independent 
Poisson and $S^\infty_+, S^\infty_-$ convergent series of independent Bernoullis.
Then for all $j\in\ZZ$ we have $\sigma_S\,\PP(S\!=\!j)\leq \eta$.
\end{corollary}

A natural question is whether such uniform bounds hold for more 
general distributions. In particular it would be interesting to characterize the distributions 
that can be obtained as limits of sums and differences of Bernoullis, beyond those 
in Corollary \ref{cor2}. In this respect we 
recall the fundamental result 
of Kintchine \cite{kin} (see also Gnedenko and Kolmogorov \cite[Theorem 2, p.115]{gnk}) 
which characterizes the limit distributions for sums of independent 
variables. The latter may or may not be Bernoullis, so that this general result provides 
only necessary conditions for our more specific question. 

\vspace{1ex}

\noindent{\sc Remark 3.} Following Remark 2, in the case of a simple Binomial 
$S\sim B(n,x)$, as well as for a single Poisson $S\sim\mathcal{P}(\lambda)$,
which is a limit of  Binomials $B(n,{\lambda\over n})$,  we have the stronger bound
$\sigma_S\,\PP(S\!=\!j)\leq{1\over\sqrt{2e}}$.

\section{Optimal bounds for fixed $n$ and $j$}
\label{s2}
Let us consider next the bound (\ref{bound}) for $n$ and $j$ fixed, namely
$$V_j^n=\max_{p\in[0,1]^n}R_j^n({}p)$$ where
$R^n_j(p)=\sigma_nP^n_j$ with $\sigma_n\!=\!\sqrt{\sum_{i=1}^np_i(1\!-\!p_i)}$ and 
$$P^n_j=\PP(S_n\!=\!j)=\sum_{|A|=j}\mbox{$\prod_{i\in A}p_i\cdot\prod_{i\not\in A}(1\!-\!p_i).$}$$
Clearly the maximum $V_j^n$ is attained and, since $R^n_j(p)$ is symmetric, any 
permutation of an optimal solution remains optimal.
Moreover, replacing each $p_i$ by $(1\!-\!p_i)$ we have the symmetry $V^n_{n-j}=V^n_j$.
It is also clear that $V^n_j$ increases with $n$ since when computing $V^{n+1}_j$ one may 
always take $p_{n+1}=0$. More generally, for $n<m$, by appropriately choosing 
$p_i\in\{0,1\}$ for $n<i\leq m$, we get
\begin{equation}\label{d1}
V^n_j\leq V_k^m\quad \forall k=j,\ldots,j+(m-n),
\end{equation}
and in particular
\begin{equation}\label{d2}
V_j^n\leq V_n^{2n}.
\end{equation}
This shows that $V^n_j$ is dominated by $V^{2n}_n$ so that 
the optimal uniform bound in (\ref{bound}) is attained as an increasing limit
$$\eta=\lim_{a\to\infty}V^{2a}_a.$$
While this was already noted in the proof of Theorem \ref{T}, the inequalities above
give a more precise picture. As a matter of fact, from (\ref{d1}) it follows that
for $n$ large, all but a small fraction of the $V^n_j$'s will be near $\eta$.
More precisely, 
\begin{proposition} \label{pp3}For each $\varepsilon>0$ there exists $n_0\in\NN$ such that
for $n\geq n_0$ and all $j$ with $\varepsilon n\leq j\leq (1-\varepsilon)n$
we have $\eta-\varepsilon\leq V^n_j\leq\eta$.
\end{proposition}
\proof Let us fix $\varepsilon>0$ and for each $n$ take $a=\lfloor\varepsilon n\rfloor$. 
For $j$ as in the statement
we have $a\leq j\leq n-a$ and we may use (\ref{d1}) to get $V^{2a}_a\leq V^n_j$.
It suffices then to choose $n$ large so that $a$ is also large enough to
ensure $V^{2a}_a>\eta-\varepsilon$.
\endproof

\noindent{\sc Remark 4.} The previous result does not hold uniformly for all $j$. In
fact, for $j=0$ one explicitly finds the optimum $p_1=\ldots=p_n=\frac{1}{2(n+1)}$
so that
$$V_0^n=V_n^n=\mbox{$\sqrt{\frac{n}{2n+1}}$}\left(1-\mbox{$\frac{1}{2(n+1)}$}\right)^{n+1}$$
which converges as $n\to \infty$ towards $\frac{1}{\sqrt{2e}}$ which is strictly smaller than $\eta$.
\\[2ex]
\noindent{\sc Conjecture 1.} {\em Numerical computations suggest that $V^n_j$ increases with $j$
for $j\leq\frac{n}{2}$ and decreases afterwards. Moreover, $V^n_j$ seems to be concave in $j$.}

\vspace{2ex}

\begin{figure}[ht]
\vspace{-1ex}
  \begin{center}
      \includegraphics[scale=0.45]{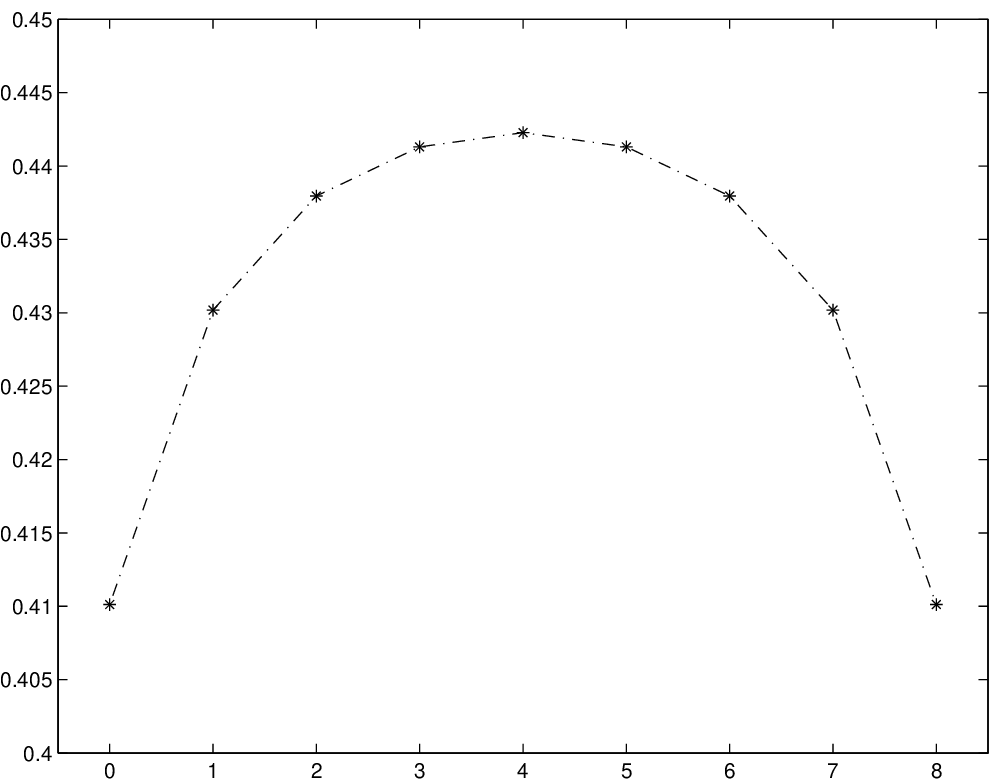}
      \vspace{-4ex}
  \end{center}
  \caption{Profile of $V^n_j$ for $j=0,\ldots,n$ (here $n=8$)\label{figura:varphi}}
\end{figure}

\subsection{Reduction to a sum of 2 Binomials}\label{dosuno}
We show next that when computing the maximum $V^n_j$ we may restrict to
$p_i$'s that take only two distinct values. This fact was 
established in \cite[Vaisman]{vai} but has not been published elsewhere.
Let us first prove that the maximum $V_j^n$ is attained with $0<p_i<1$ 
for  all $i=1,\ldots,n$, which 
is equivalent to showing that  the inequalities (\ref{d1}) are strict. We  exploit the following properties.

\begin{lemma} \label{lema3} Let $p$ be optimal for $V^n_j$ and let $P^n_k=\PP(S_n=k)$ 
be the corresponding distribution with mean $\mu_n$ and variance
$\sigma_n^2$.
If $0<p_i<1$ for $i=1,\ldots,n$ then
\begin{itemize}
\item[(a)] $2\sigma_n^2[P^n_j-P^n_{j+1}]=\frac{\mu_n}{j+1}P^n_j$,
\item[(b)] $2\sigma_n^2[P^n_j-P^n_{j-1}]=\frac{n-\mu_n}{n-j+1}P^n_j$,
\item[({}c)] $S_n$ has a unique mode at $j$,
\item[(d)] $\sigma_n^2[j-\mu_n]+\sum_{i=1}^np_i(1\!-\!p_i)(\mbox{$\frac{1}{2}$}-p_i)=0$,
\item[(e)] $j-\frac{1}{2}<\mu_n<j+\frac{1}{2}$.
\end{itemize}
\end{lemma}
\proof
Letting $S^i=\sum_{j\neq i}X_j$ and $Q_k^i=\PP(S^i\!=\!k)$ we have the known identities
\begin{eqnarray}
\label{ll1}  P^n_k&=&p_iQ^i_{k-1}+(1\!-\!p_i)Q_k^i,\\
\label{ll0}  kP^n_k&=&\mbox{$\sum_{i=1}^np_iQ^i_{k-1}.$}
\end{eqnarray}
Since $p$ is an interior optimal solution, the optimality conditions give
\begin{equation}\label{ll12}
\frac{\partial R_j^n}{\partial p_i}=\frac{1-2p_i}{2\sigma_n}P^n_j+\sigma_n[Q^i_{j-1}-Q^i_j]=0
\end{equation}
which multiplied by $p_i/\sigma_n$ and summed over $i$ yields
$$\mbox{$\frac{1}{2\sigma_n^2}$}P^n_j\sum_{i=1}^np_i(1\!-\!2p_i)+\sum_{i=1}^np_i[Q^i_{j-1}-Q^i_j]=0.$$
The first sum is just $2\sigma_n^2-\mu_n$, while (\ref{ll0}) shows that the 
second sum is equal to $jP^n_j-(j+1)P^n_{j+1}$, so that rearranging terms we get {\em (a)}.

Property {\em (b)} follows from {\em (a)} applied to $\bar p_i=1\!-\!p_i$ which is optimal for $V^n_{n-j}$,
while {\em ({}c)} is a consequence of the fact that the distribution of $S_n$ is unimodal
combined with {\em (a)} and {\em (b)} that give respectively $P^n_j>P_{j+1}^n$ and $P^n_j>P^n_{j-1}$.

To prove {\em (d)} we multiply  (\ref{ll12}) by $\sigma_n\,p_i(1\!-\!p_i)$ and sum over $i$ to get
\begin{equation}\label{ll2}
P^n_j\sum_{i=1}^np_i(1\!-\!p_i)(\mbox{$\frac{1}{2}$}-p_i)+\sigma_n^2\sum_{i=1}^np_i(1\!-\!p_i)[Q^i_{j-1}-Q^i_j]=0.
\end{equation}
Now, from (\ref{ll1}) we have $P^n_j-Q^i_{j-1}=(1\!-\!p_i)[Q^i_j-Q^i_{j-1}]$ and using (\ref{ll0}) we get
$$\sum_{i=1}^np_i(1\!-\!p_i)[Q^i_{j-1}-Q^i_j]=\sum_{i=1}^np_i[Q^i_{j-1}-P^n_j]=[\mbox{$j-\sum_{i=1}^np_i$}]P^n_j$$
which plugged into (\ref{ll2}) and simplifying by $P^n_j$ yields {\em ({}d)}.

Finally, to prove {\em ({}e)} it suffices to observe that {\em ({}d)} can be rewritten as
\begin{eqnarray*}
\sigma_n^2[\mbox{$j+\frac{1}{2}-\mu_n$}]&=&\hspace{1.7ex}\mbox{$\sum_{i=1}^np_i^2(1\!-\!p_i)$}~>0,\\
\sigma_n^2[\mbox{$j-\frac{1}{2}-\mu_n$}]&=&-\mbox{$\sum_{i=1}^np_i(1\!-\!p_i)^2$}<0.
\end{eqnarray*}

\vspace{-5ex}
\endproof

\begin{proposition} If $p$ attains the maximum $V_j^n$ then $0<p_i<1$ for $i=1,\ldots,n$.
Moreover, for $n<m$ all the inequalities {\em (\ref{d1})} are strict.
\end{proposition}
\proof The strict inequality in (\ref{d1}) is a
direct consequence of the fact that the optimal $p_i$'s do not take the values 0 nor 1.
We prove the latter by induction in $n$. The property clearly holds for $n=1$.
Assume that it holds for a given $n$ and let us prove 
it for $n+1$. Take $p$ optimal for $V^{n+1}_j$ and suppose for a contradiction
that it has a null component, say $p_{n+1}=0$. In this case 
$V^{n+1}_j=V^n_j$ and $(p_1,\ldots,p_n)$ is optimal for $V_j^n$
so the induction hypothesis yields $0<p_i<1$ for $i=1,\ldots,n$. Denoting as before
$P^n_k=\PP(S_n=k)$ and using properties {\em (b)} and {\em (e)} in Lemma \ref{lema3} we obtain
$$\frac{\partial R^{n+1}_j}{\partial p_{n+1}}=\frac{1}{2\sigma_n}P^n_j+\sigma_n[P^n_{j-1}-P^n_j]=\frac{P^n_j}{2\sigma_n}\frac{\mu_n-j+1}{n-j+1}>0$$
which shows that $p_{n+1}=0$ cannot be a maximizer. 
This same argument applied to $\bar p_i=1\!-\!p_i$, which is optimal for $V^{n+1}_{n-j}$,
shows that no $p_i$ can be equal to 1.
\endproof

With this result we may now show that when computing 
$V^n_j$ one may restrict to a sum of two Binomials.
\begin{proposition} There exists $p$ optimal for $V^n_j$ which takes
at most two distinct values $p_i\in\{\alpha,\beta\}$ with $\alpha,\beta\in (0,1)$.
In other words, the maximum is attained for $S_n=U+V$ with $U\sim B(a,\alpha)$ 
and $V\sim B(b,\beta)$ independent Binomials with $a+b=n$. More explicitly,
denoting $b^n_k(x)={n\choose k}x^k(1\!-\!x)^{n-k}$, we have
\begin{equation}\label{optimal}
V^n_j=\max_{\stackrel{\mbox{\scriptsize $ a\!+\!b=n$}}{\alpha,\beta\in (0,1)}}\sqrt{a\,\alpha(1\!-\!\alpha)+b\,\beta(1\!-\!\beta)}\;\sum_{k=0}^jb^a_k(\alpha)b^b_{j-k}(\beta).
\end{equation}

\end{proposition}

\proof Take $p$ an optimal solution for $V^n_j$ with the smallest product $\Pi_{i=1}^np_i$.
We claim that this $p$ takes at most two values.
Assume by contradiction that it has 3 different entries $0\!<\!p_r\!<\!p_s\!<\!p_t\!<\!1$.
Denoting $\tilde p=(p_i)_{i\neq r,s,t}$ and $\bar p_i=1-p_i$ we have
\begin{eqnarray*}
P^n_j&=&p_rp_sp_tP^{n-3}_{j-3}(\tilde p)+[\bar p_rp_sp_t+p_r\bar p_sp_t+p_rp_s\bar p_t]P^{n-3}_{j-2}(\tilde p)\\
&&\hspace{-0.5cm}\mbox{}+[\bar p_r\bar p_sp_t+\bar p_rp_s\bar p_t+p_r\bar p_s\bar p_t]P^{n-3}_{j-1}(\tilde p)+\bar p_r\bar p_s\bar p_tP^{n-3}_j(\tilde p)
\end{eqnarray*}
which may be rewritten as $P^n_j=F(p_r,p_s,p_t)$ where
$$F(x,y,z)=Axyz+B(xy\!+\!xz\!+\!yz)+C(x\!+\!y\!+\!z)+D$$
with coefficients $A,B,C,D$ that depend only on $\tilde p$. 
Setting $\kappa=P^n_j/2\sigma_n^2$ the optimality conditions for $p$ yield
\begin{eqnarray*}
\kappa(1\!-\!2p_r) + Ap_sp_t+B(p_s+p_t)+C&=& 0\\
 \kappa(1\!-\!2p_s)+Ap_rp_t+B(p_r+p_t)+C&=&0\\
\kappa(1\!-\!2p_t)+Ap_rp_s+B(p_r+p_s)+C&=& 0.
\end{eqnarray*}
Substracting the first two equations and simplifying by $p_s-p_r\neq 0$
we get $2 \kappa +Ap_t+B=0$, while the second and third
equations give $2 \kappa+Ap_r+B=0 $. Hence $Ap_t=Ap_r$ and since $p_r\neq p_t$ 
we conclude $A=0$ so that
$F(x,y,z)$ depends only on $x+y+z$ and $xy+xz+yz$. Moreover, since we also have
$$x(1\!-\!x)+y(1\!-\!y)+z(1\!-\!z)=(x\!+\!y\!+\!z)-(x\!+\!y\!+\!z)^2+2(xy\!+\!xz\!+\!yz)$$ 
it follows that $R^n_j(x,y,z,\tilde p)$ is constant over the set defined by the equations
$x\!+\!y\!+\!z=p_r\!+\!p_s\!+\!p_t$ and $xy\!+\!xz\!+\!yz=p_rp_s\!+\!p_rp_t\!+\!p_sp_t$.
Thus, any such vector $(x,y,z,\tilde p)$ maximizes $R^n_j(\cdot)$ and our choice of $p$ 
implies that $(p_r,p_s,p_t)$ solves
$$\left\{\begin{array}{l}
\min~~xyz\\
s.t.~~(x,y,z) \in [0,1]^3\\
\begin{array}{lcl}
       x+y+z&=&p_r+p_s+p_t\\
       xy+xz+yz&=&p_rp_s+p_rp_t+p_sp_t.
      \end{array} 
\end{array}\right.
$$
Since $p_r,p_s,p_t$ are different, the gradients of the two equality constraints at
this optimal point are linearly independent, while the inequality constraints are non-binding.
Hence the Mangasarian-Fromovitz constraint qualification holds and we may
find Lagrange multipliers $\lambda$ and $\mu$ such that
\begin{eqnarray*}
p_sp_t&=& \lambda +\mu(p_s+p_t)\\
p_rp_t&=& \lambda + \mu(p_r+p_t)\\
p_rp_s&=& \lambda + \mu(p_r+p_s).
\end{eqnarray*}
Substracting the first two equations and simplifying by $p_s\!-p_r\!\neq\! 0$ we get $p_t\!=\!\mu$,
and similarly $p_r=\mu$ and $p_s=\mu$.
This contradiction shows that $A$ 
cannot be 0, and therefore the assumption $0<p_r<p_s<p_t<1$ was absurd.  
\endproof

\vspace{1ex}
\noindent
 {\sc Conjecture 2:} {\em We conjecture that for each $n$ and $j$ there exist 
 unique values $\alpha <\frac{1}{2}<\beta $ such that the maximizers of  $V^n_j$ 
 are precisely the vectors with exactly $(n-j)$ components equal to $\alpha$ and $j$ 
components equal to $\beta$.}

\vspace{2ex}
Remark 4 shows that this property holds for $V^n_0$ and $V^n_n$.
In the next section we prove that it also holds for the dominating values 
$V^{2a}_a$ and,  moreover, in this case $\alpha+\beta=1$.

\subsection{Computation of the sharp bound $V^{2a}_a$}
Proceeding as in the proof of Theorem \ref{T}, setting $z_i=p_i(1\!-\!p_i)\in [0,
\frac{1}{4}]$ we have
$$\PP(S_n=j)\leq\frac{1}{2\pi}\int_0^{2\pi}\prod_{i=1}^n\sqrt{1+2z_i(\cos\theta-1)}\;d\theta$$
so that using the change of variables $\theta=2\xi$ we get
\begin{equation}\label{cr1}
R^n_j({}p)\leq\Phi_n(z):=\sqrt{\mbox{$\sum_{i=1}^nz_i$}}\int_0^{\pi}\!\!\frac{1}{\pi}\prod_{i=1}^n\sqrt{1\!-\!4z_i\sin^2\!\xi}\;d\xi.
\end{equation}
Therefore, denoting $\eta_n$ the maximum of $\Phi_n(z)$ for $z\in[0,\frac{1}{4}]^n$ we have
$V^n_j\leq\eta_n$. We will prove that for $n=2a$ and $j=a$ we have the equality $V^{2a}_a=\eta_{2a}$.
Moreover, we will show that $V^{2a}_a$ has a unique maximizer $p$ up to permutation.
We begin by characterizing the maximum of $\Phi_n$ for $n=2a$. 
\begin{lemma}\label{l7}
For $n=2a$ the map $\Phi_n(z)$ has a unique maximizer  $z\in[0,\frac{1}{4}]^n$
which is of the form $z_i=\bar z$ for all $i=1,\ldots,n$ with $\bar z\in (0,\frac{1}{4})$.
\end{lemma}
\proof 
Let $z\in[0,\frac{1}{4}]^n$ maximize $\Phi_n$. Denoting $Q_z(\xi)=\frac{1}{\pi}\prod_{i=1}^n\sqrt{1\!-\!4z_i\sin^2\!\xi}$
and $s=\sum_{i=1}^nz_i$ we have
\begin{equation}\label{derstat}
\frac{\partial\Phi_n}{\partial z_k}(z)=\frac{1}{2\sqrt{s}}\int_0^{\pi}\!\!\!Q_z(\xi)d\xi-2\sqrt{s}\int_0^{\pi}\!\!\!Q_z(\xi)\frac{\sin^2\!\xi}{1\!-\!4z_k\sin^2\!\xi}d\xi.
\end{equation}
Clearly the $z_i$'s are not all 0. They cannot be all equal to $\frac{1}{4}$ either, since
in that case we would have $Q_z(\xi)=\frac{1}{\pi}\cos^{2a}\!\xi$ so that letting
$$J^{k,i}=\frac{1}{\pi}\int_0^\pi\cos^{2k}\!\xi\sin^{2i}\!\xi\,d\xi=\frac{1}{4^{k+i}}\frac{(2k)!(2i)!}{k!i!(k+i)!}$$
we get the following inequality that contradicts optimality
$$\frac{\partial\Phi_n}{\partial z_k}(z)=\mbox{$\frac{1}{\sqrt{2a}}$}\,J^{a,0}-\sqrt{2a}\,J^{a-1,1}=-\mbox{$\sqrt{2a}\,\frac{(2a-2)!}{(2^a a!)^2}$}<0.
$$
Hence, there is some component $z_k\in (0,\frac{1}{4})$ for which we have
$\frac{\partial\Phi_n}{\partial z_k}=0$. Using (\ref{derstat}) we then get
$\frac{\partial\Phi_n}{\partial z_i}=\frac{\partial\Phi_n}{\partial z_i}-\frac{\partial\Phi_n}{\partial z_k}=(z_k-z_i)Q_{ik}$ with
$$Q_{ik}=
\int_0^{\pi}\!\!\!\frac{8\sqrt{s}\,Q_z(\xi)\sin^4\!\xi}{[1\!-\!4z_k\sin^2\!\xi][1\!-\!4z_i\sin^2\!\xi]}d\xi.$$
Note that if $z_i=\frac{1}{4}$ and all the other components are smaller than $\frac{1}{4}$ this integral
is singular and $Q_{ik}=\infty$. In all other cases $Q_{ik}$ is finite and strictly positive.
Hence, if $z_i=0$ we get $\frac{\partial\Phi_n}{\partial z_i}>0$ contradicting optimality, while $z_i=\frac{1}{4}$ 
leads to the contradiction $\frac{\partial\Phi_n}{\partial z_i}<0$.
Therefore $z_i\in(0,\frac{1}{4})$ so that $\frac{\partial\Phi_n}{\partial z_i}= 0$,
which implies $z_i=z_k$ and then $z=(\bar z,\ldots,\bar z)$ with $\bar z\in (0,\frac{1}{4})$. 

It remains to show that $\bar z$ is unique. A routine calculation shows that the optimality condition 
$\frac{\partial\Phi_n}{\partial z_k}(\bar z,\ldots,\bar z)=0$ is equivalent to $h(\bar z)=0$ where $h(z)$ is the polynomial
$$h(z)=\frac{2a\!+\!1}{\pi}\int_0^{\pi}\!\![1\!-\!4z\sin^2\!\xi]^a d\xi-\frac{2a}{\pi}\int_0^{\pi}\!\![1\!-\!4z\sin^2\!\xi]^{a-1} d\xi.$$
According to Fourier-Budan's Theorem the number of roots in $(0,\frac{1}{4})$ 
is at most $\sigma(0)-\sigma(\frac{1}{4})$ where $\sigma(z)$ denotes the number of sign changes
in the sequence $h(z),h'(z),h''(z),\ldots,h^{(a)}(z)$. By direct calculation we get 
$$h^{(i)}(z)=\frac{(-4)^ia!}{(a-i)!}\frac{1}{\pi}\int_0^\pi[1-4z\sin^2\!\xi]^{a-i}\sin^{2i}\!\xi\left[2a+1-\frac{2(a-i)}{1\!-\!4z\sin^2\!\xi}\right]d\xi.$$
Hence $h^{(i)}(0)=\frac{(-4)^ia!}{(a-i)!}(2i+1)J^{0,i}$ which alternates sign so that $\sigma(0)=a$.
Also 
\begin{eqnarray*}
h^{(i)}(\mbox{$\frac{1}{4}$})&=&\frac{(-4)^ia!}{(a-i)!}\frac{1}{\pi}\int_0^\pi\cos^{2(a-i)}\!\xi\sin^{2i}\!\xi\left[2a+1-\frac{2(a-i)}{\cos^2\!\xi}\right]d\xi\\
&=&\frac{(-4)^ia!}{(a-i)!}\left[(2a\!+\!1)J^{a-i,i}-2(a\!-\!i)J^{a-i-1,i}\right]
\end{eqnarray*}
where for $i=a$ the last term $J^{-1,a}$ is interpreted as 0. Explicitely we get
$$h^{(i)}(\mbox{$\frac{1}{4}$})=
\left\{
\begin{array}{ll}
(-1)^{i+1}4^{i-a}\frac{(2a-2i)!(2i)!}{(a-i)!^2i!}\frac{2i+1}{2a-2i-1}&\mbox{ for $i<a$}\\[1ex]
(-1)^a(2a\!+\!1)\frac{(2a)!}{a!}&\mbox{ for $i=a$}
\end{array}
\right.$$ 
which yields $\sigma(\frac{1}{4})=a-1$. Hence $\sigma(0)-\sigma(\frac{1}{4})=1$
and therefore $h$ has exactly one root $\bar z$ in $(0,\frac{1}{4})$.
\endproof

According to Lemma \ref{l7}, the maximum of $\Phi_n$ for $n=2a$ is given by
\begin{equation}\label{cota2a}
\eta_{2a}=\max_{0\leq z\leq\frac{1}{4}}\frac{\sqrt{2az}}{\pi}\int_0^\pi(1\!-\!4z\sin^2\!\xi)^a\,d\xi
\end{equation}
which is attained at a unique $\bar z\in (0,\frac{1}{4})$. Note that using the change of variables 
$t=\sin^2\!\xi$ the integral can be expressed using the hypergeometric function $_2F_1$, 
and the function to be maximized is $\sqrt{2az}\,{}_2F_1(-a,\frac{1}{2};1;4z)$.
\begin{theorem} For $n=2a$ and $j=a$ we have $V^{2a}_a=\eta_{2a}$. Moreover, let $\alpha$ be the unique solution of 
$x(1\!-\!x)=\bar z$ with $\alpha\in(0,\frac{1}{2})$ and $\bar z$ the 
solution of {\em (\ref{cota2a})}. Then the maximizers of $V^{2a}_a$ are precisely the vectors 
$p$ with half of its components equal to $\alpha$ and the other half equal to $1\!-\!\alpha$.
\end{theorem}
\proof From (\ref{cr1}) we have $V^{2a}_a\leq\eta_{2a}$, while taking
$p$ with half of the $p_i$'s equal to $\alpha$ and the other half $1\!-\!\alpha$ equation (\ref{pgf}) 
gives $R^{2a}_a({}p)=\Phi_n(\bar z,\ldots,\bar z)=\eta_{2a}$
so that in fact $V^{2a}_a=\eta_{2a}$ and this $p$ is optimal for $V^{2a}_a$.

Now take any maximizer $p$ of $V^{2a}_a$ and define $z\in[0,\frac{1}{4}]^n$ 
by $z_i=p_i(1\!-\!p_i)$. Then $V^{2a}_a=R^{2a}_a({}p)\leq\Phi_n(z)\leq\eta_{2a}$
so that $z$ is optimal for $\Phi_n$ and therefore $z_i=\bar z$ 
for all $i=1,\ldots,n$. It follows that $p_i(1\!-\!p_i)=\bar z$ and then all the components $p_i$ are either $\alpha$ 
or $1\!-\!\alpha$. From Lemma \ref{lema3}{\em (d)} it readily follows that there are exactly 
half of them of each type.
\endproof

\vspace{1ex}
We conclude this section with a quantitative version of Proposition \ref{pp3}.
\begin{proposition} For each  $n\geq 1$ and $0<j<n$ let $k=\min(j,n-j)$. Then
\begin{equation}\label{final2}
\mbox{$(1\!-\!\frac{1}{2k})$} \eta< V^n_j<\eta.
\end{equation}
In particular $V^{2a}_a$ converges to $\eta$ at rate $O(\frac{1}{a})$.
\end{proposition}
\proof Let $x_0\sim 0.78997786$ be the point where the maximum $\eta=\sqrt{x_0}\,e^{-x_0}I_0(x_0)$
is attained.
Taking  $z=x_0/2a$ in (\ref{cota2a}) 
and using the inequality $(1\!-\!\frac{y}{a})^{a-1}> e^{-y}$ 
which holds for all $a\geq 1$ and $y\in(0,2x_0)$, it follows that
$$\eta_{2a}\geq\frac{\sqrt{x_0}}{\pi}\!\int_0^\pi\!\!\mbox{$\left(1\!-\!\frac{2x_0\sin^2\!\xi}{a}\right)^a\,d\xi$}
>\frac{\sqrt{x_0}}{\pi}\!\!\int_0^\pi\!\!\!\mbox{$(1-\frac{2x_0\sin^2\!\xi}{a})e^{-2x_0\sin^2\!\xi}\,d\xi$}.$$
Using the change of variables $\theta=2\xi$, the latter can be expressed in terms of the
modified Bessel functions $I_0(\cdot)$ and $I_1(\cdot)$ as
\begin{eqnarray*}
\eta_{2a}&>&\frac{\sqrt{x_0}}{2\pi}\!\!\int_0^{2\pi}\!\!\!\mbox{$(1+\frac{x_0}{a}(\cos\theta\!-\!1))e^{x_0(\cos\theta-1)}\,d\theta$}\\
&=&\sqrt{x_0}\,e^{-x_0}\left[ (1-\mbox{$\frac{x_0}{a}$})I_0(x_0)+\mbox{$\frac{x_0}{a}$}I_1(x_0)\right].
\end{eqnarray*}
Now, $I_1(\cdot)=I_0'(\cdot)$ while by optimality the derivative of $\sqrt{x}\,e^{-x}I_0(x)$ 
vanishes at $x_0$ which yields $x_0I_1(x_0)=x_0I_0'(x_0)=(x_0-\frac{1}{2})I_0(x_0)$, and therefore we get
\begin{equation}\label{final1}
V^{2a}_a=\eta_{2a}>\mbox{$(1-\frac{1}{2a})$}\sqrt{x_0}\,e^{-x_0}I_0(x_0)=\mbox{$(1-\frac{1}{2a})$}\eta.
\end{equation}

From this inequality we easily derive (\ref{final2}). Indeed, if $j\leq n-j$ then $2j\leq n$ and (\ref{d1}) implies $V^{2j}_j\leq V^n_j$ 
which combined with 
(\ref{final1}) yields $(1\!-\!\frac{1}{2j})\eta<V^n_j$. The case $j\geq n-j$ follows from this 
by exploiting the symmetry $V^n_j=V^n_{n-j}$. 
\endproof

\section{An application to fixed point iterations}\label{appl}
Let us illustrate how Theorem \ref{T} can be used to study the rate of 
convergence of fixed point iterations \cite{kra, man}. Namely, let $(E,\|\cdot\|)$
 be a normed vector space and $T:E\to E$ a non-expansive map, that is, 
 $\|T(x)-T(y)\|\leq\|x-y\|$ for all $x,y\in E$, with a nonempty set of fixed points $\mbox{Fix}(T)$. 
Consider the Krasnosel'ski\v{\i}-Mann iteration
  $$x_n=(1-\alpha_n)x_{n-1}+\alpha_n Tx_{n-1}$$
with $x_0\in E$ given and $0<\alpha_n<1$.
In  \cite{bb}, Baillon and Bruck conjectured the existence of a universal constant $C$
such that 
\begin{equation}\label{bnd}
\|x_n-Tx_n\|\leq C\frac{\mbox{dist}(x_0,\mbox{Fix}(T))}{\sqrt{\sum_{i=1}^n\alpha_i(1-\alpha_i)}}
\end{equation}
proving this bound with $C=2/\sqrt{\pi}\sim 1.1284$ for $\alpha_i\equiv\alpha$ constant.
The general case with non-constant $\alpha_i$'s was recently settled in \cite{cc}
with this same $C$, while \cite[Vaisman]{vai} proved that it
holds with $C=1$ when $E$ is a Hilbert space. Here we use Theorem \ref{T} to find a slightly 
improved bound for affine maps in general normed spaces.

 \begin{proposition}\label{p11} Let  $T(x)\!=\!a\!+\!Lx$ with $L\!:\!E\!\to\! E$ linear and non-expansive. Then 
 {\em (\ref{bnd})} holds with $C=2\eta\sim  0.9376$ where $\eta$ is given
 by {\em (\ref{bnd1})}.
 \end{proposition}
 \proof A simple inductive argument shows that 
 $x_n=\sum_{j=0}^n\pi_j^n\,T^jx_0$
 where the coefficients $\pi_j^n$ satisfy the recursion
 $\pi_j^n=\alpha_n\pi_{j-1}^{n-1}+(1-\alpha_n)\pi_j^{n-1}$.
 Notice that $\pi_j^n=\PP(S_n=j)$ where $S_n\!=\!X_1+\cdots+X_n$ is a sum of independent 
Bernoullis with $\PP(X_i\!=\!1)=\alpha_i$. In particular $\sum_{j=0}^n\pi_j^n=1$
so that for each $y\in\mbox{ Fix}(T)$ we have $x_n-y=\sum_{j=0}^n\pi_j^nT^j(x_0-y)$,
and since $x_n-Tx_n=(x_n-y)- T(x_n-y)$ the triangle inequality implies
$$\|x_n-Tx_n\|\leq \sum_{j=0}^{n+1}|\pi_j^n-\pi_{j-1}^n|\|T^j(x_0-y)\|\leq \|x_0-y\| \sum_{j=0}^{n+1}|\pi_j^n-\pi_{j-1}^n|.$$
Since the distribution of $S_n$ is unimodal the latter sum
is $2\max_j\pi_j^n$. The conclusion follows by using (\ref{bound}) and taking infimum
over $y\in \mbox{Fix}(T)$.
 \endproof

 \vspace{4ex}
 \noindent{\bf Acknowledgements:} We are indebted to Professor David McDonald for 
 helpful discussions on the connection of our main result with the local limit theorem,
 as well as to Professor  Michel Weber for pointing out the Kolmogorov-Rogozin inequality.
 We also thank an anonymous referee for very useful suggestions that contributed to improve 
 the presentation. Roberto Cominetti was supported by FONDECYT Grant 1100046 (CONICYT-Chile),
 as well as Nucleo Milenio Informaci\'on y Coordinaci\'on en Redes ICM/FIC P10-024F.
This work was completed during a visit of Jean-Bernard Baillon to Universidad de Chile,
which was supported by FONDECYT Grant 1130564.

\end{document}